\documentclass[11pt, onecolumn, a4paper]{article}
\usepackage{amsmath}
\usepackage{amssymb}
\usepackage{amsthm}
\usepackage{lipsum}
\usepackage{mathptmx}

\usepackage{hyperref} 
\usepackage{url} 
\usepackage{bm}

\usepackage[rm,up,sc,topmarks,calcwidth,pagestyles]{titlesec}
\titleformat{\section}{\Large\bfseries\rmfamily}{\thesection}{1em}{}
\titleformat{\subsection}{\large\bfseries\rmfamily}{\thesubsection}{1em}{}
\titleformat{\subsubsection}{\large\it\rmfamily}{\thesubsubsection}{1em}{}

\theoremstyle{theorem}
\theoremstyle{proposition}
\theoremstyle{definition}

\newtheorem{proposition}{Proposition}

\usepackage{cite}
\usepackage[dvipdfmx]{graphicx}
\usepackage{textcomp}
\usepackage{url}

\usepackage{color}

\usepackage{algorithm}
\usepackage{algorithmic}

\title{Ranks of Checkered Pattern Matrices with Applications to Information Embedding and Retrieving}
\author{Hideo Hirose\footnote{Chuo University \& Kurume University, Japan} }
\date{}
\begin{document}
\maketitle
\thispagestyle{empty}

\section*{Abstract}
Checkered patterns are characterized by their square structure and the use of only two distinct colors. These colors are typically represented by two types of numerical sets: $\{1,0\}$ and $\{1,-1\}$. Matrices based on $\{1,0\}$ may seem identical to those based on $\{1,-1\}$ when forming checkered patterns because the only difference is that the numbers $0$ are changed to $-1$. However, these two kinds of matrices are completely different in a mathematical sense because a matrix using $\{1, 0\}$ has a rank of $2$ and a matrix using $\{1, -1\}$ has a rank of $1$. Knowing this difference in advance allows us to reduce the computational effort required for matrix operations such as information embedding and retrieving.
\\[2mm]
{\it Keywords: } checkered pattern, matrix rank, bijective transformation, matrix decomposition, information embedding and retrieving.

\section{Introduction}

Checkered patterns are found everywhere not only in human-made objects, including textile designs, checkerboards in the board game chess, and flags in auto racing, but also in nature~\cite{ConnerCollinsSimberloff2013}. 
Although 
the checkered pattern, widely used in architecture and interior fixtures, embodies a deceptively simple design, its bold contrast and rhythmic repetition lend it both familiarity and timeless appeal, thereby making it highly significant and appealing in aesthetic contexts.
For example, 
Figure 1 shows a traditional Japanese {\it fusuma} painting featuring the checkered pattern, as well as a designer hotel corridor that incorporates the same motif.
Even with only two distinctive colors (bright color and dark color), they convey such visual force and aesthetic appeal that they leave a lasting impression on the observer. 
 
\begin{figure}[htb]
\begin{center}
\includegraphics[scale=0.4]{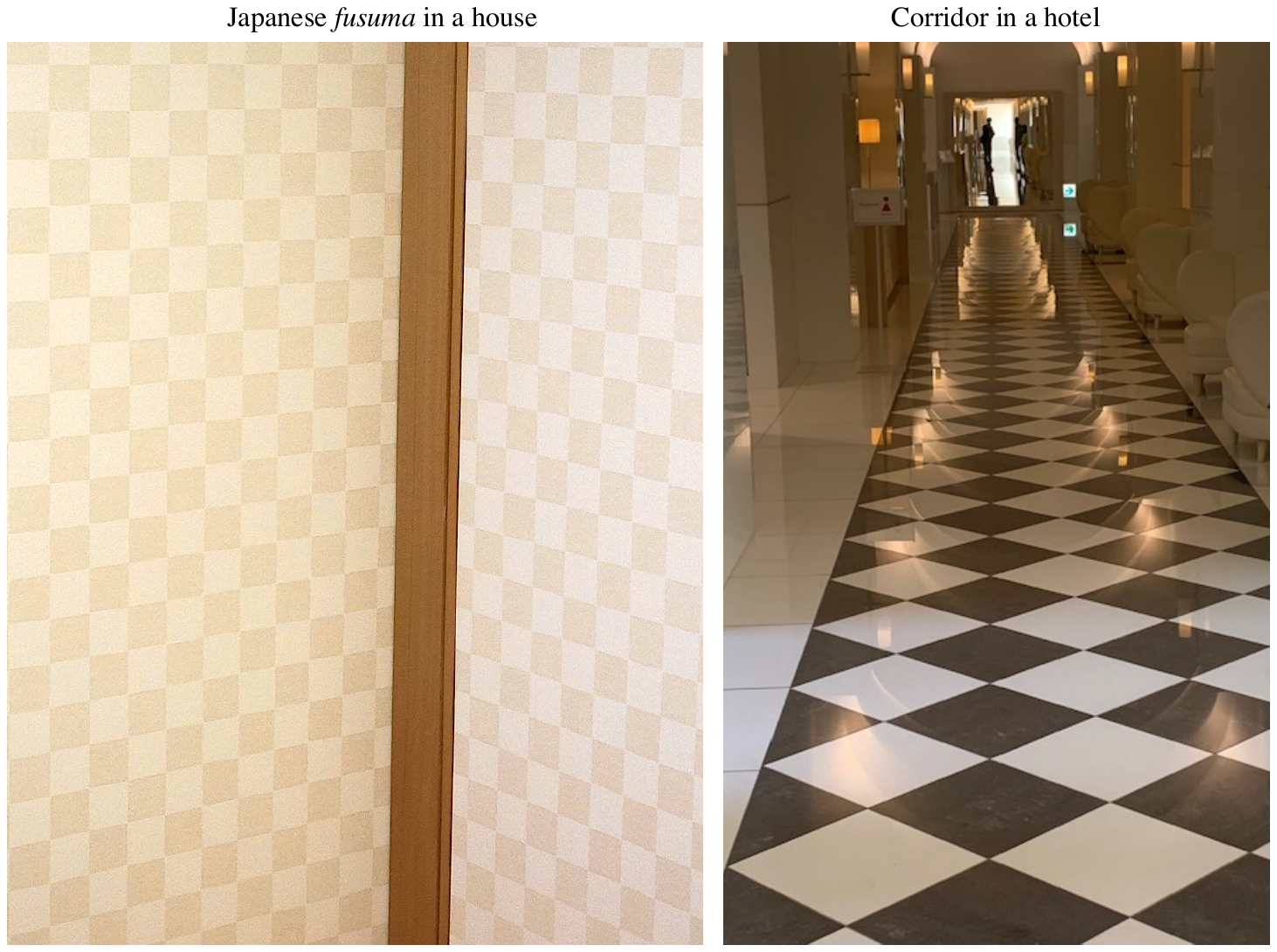}
\end{center}
\caption{Checkered patterns seen in a Japanese traditional {\it fusuma} painting and a designer hotel corridor.} 
\label{fig:FusumaCorridor}
\end{figure}

In mathematical sciences, checkered patterns have found applications in numerical computations~\cite{DelgadoOreraPena2024, CrasterObnosov2001} and in graph~\cite{SilverWilliams2019}. Beyond mathematics, they also appear in the social and ecological sciences, where they model patterns of social interaction~\cite{Sakoda1971, Diamond1975, StoneRoberts1990}. In technology, their utility extends further still, contributing to robust camera calibration~\cite{YuPeng2006, WangZhengZhang2021}, the design of efficient solar cells~\cite{Donaldson2021}, and techniques in steganography and digital watermarking~\cite{ItoKanazawa2023}.

What checkered patterns share in common is that they are square in form and encode only two distinct values.
These binary values are commonly expressed as either $\{1,0\}$ or $\{1,-1\}$. The former typically appears in digital data (bit representation), probability theory (success or failure), and logic (true or false), whereas the latter is more often employed in signal processing and machine learning. Because there is a one-to-one correspondence between $\{1,0\}$ and $\{1,-1\}$, the two representations are mathematically equivalent.

Using these two sets, we construct two types of checkered pattern $m \times n$ matrices $A$ and $B$,
\begin{eqnarray}
\label{eq:intro_1}
   A = 
   \begin{pmatrix}
        1 & 0 & 1 & 0 &\cdots  \\
        0 & 1 & 0 & 1 &  \cdots  \\
        1 & 0 & 1 & 0 & \cdots  \\
        0 & 1 & 0 & 1 &  \cdots  \\
        \vdots & \vdots & \vdots & \vdots & \ddots \\
   \end{pmatrix}, \ 
   B = 
   \begin{pmatrix}
        1 & -1 & 1 & -1 & \cdots  \\
        -1 & 1 & -1 & 1 & \cdots  \\
        1 & -1 & 1 & -1 & \cdots  \\
        -1 & 1 & -1 & 1 & \cdots  \\
        \vdots & \vdots & \vdots & \vdots & \ddots \\
   \end{pmatrix}.    
\end{eqnarray}
These two patterns 
may seem identical to each other
because the two forms are connected by a linear-affine transformation such that 
\begin{eqnarray}
\label{eq:intro_2}
A = \frac{1}{2}(B + C)
\end{eqnarray}
where $C$ is the all-ones $m \times n$ matrix expressed as
\begin{eqnarray}
\label{eq:sec1_3}
   C = 
   \begin{pmatrix}
        1 & 1 & 1 & 1 &\cdots  \\
        1 & 1 & 1 & 1 &\cdots  \\
        1 & 1 & 1 & 1 &\cdots  \\
        1 & 1 & 1 & 1 &\cdots  \\
        \vdots & \vdots & \vdots & \vdots & \ddots \\
   \end{pmatrix}.    
\end{eqnarray}
That is, both transformations ($A$ and $B$) are bijective linear maps up to a constant shift.
Therefore, we are tempted to conclude that both matrices possess the same characteristics
such that the transformation preserves the rank, and null space.
Nevertheless, in fact, they are completely different in a mathematical sense because matrix $A$ has a rank of $2$ and matrix $B$ has a rank of $1$. 
This is a kind of paradox.
Surprisingly, despite its simplicity, this apparent contradiction has received no attention in the literature to date, as far as I know. 
This seemingly limiting rank-reduction property can actually be harnessed to reduce computational cost in matrix operations.

\section{Rank of Checkered Pattern Matrix}

\begin{proposition}[Rank of Matrix]
The rank of matrix $A$ is $2$, and the rank of matrix $B$ is $1$.
\end{proposition}

\begin{proof}
Applying row operations to $A$, $A$ becomes the following matrix
\begin{eqnarray}
\label{eq:sec1_1}
   \begin{pmatrix}
        1 & 0 & 1 & 0 & \cdots  \\
        0 & 1 & 0 & 1 & \cdots  \\
        0 & 0 & 0 & 0 &\cdots  \\
        0 & 0 & 0 & 0 &\cdots  \\
        \vdots & \vdots & \vdots & \vdots & \ddots \\
   \end{pmatrix}    
\end{eqnarray}
of rank $2$.

Similarly, using row operations to $B$, $B$ becomes the following matrix
\begin{eqnarray}
\label{eq:sec1_2}
   \begin{pmatrix}  
        1 & -1 & 1 & -1 & \cdots  \\
        0 & 0 & 0 & 0 &\cdots  \\
        0 & 0 & 0 & 0 &\cdots  \\
        0 & 0 & 0 & 0 &\cdots  \\
        \vdots & \vdots & \vdots & \vdots & \ddots \\
   \end{pmatrix}    
\end{eqnarray}
of rank $1$.
Such a proof follows essentially the $LU$ decomposition procedure.
\end{proof}

\medskip
To find factors that lead to the difference in rank between the two matrices,
we extend equation (\ref{eq:intro_2}) into a more general form by means of a linear combination parameter $-\infty < \alpha < \infty$. A linearly combined matrix using $B$ and $C$ is
\begin{eqnarray} 
\label{eq:sec1_4}
   (1 - \alpha) B + \alpha C 
   &=& 
   \begin{pmatrix}
        1  & -1 + 2\alpha & 1 & -1 + 2\alpha & \cdots  \\
        -1 + 2\alpha & 1 & -1 + 2\alpha & 1 & \cdots  \\
        1 & -1 + 2\alpha & 1 & -1 + 2\alpha & \cdots  \\
        -1 + 2\alpha & 1 & -1 + 2\alpha & 1 & \cdots  \\
        \vdots & \vdots & \vdots & \vdots & \ddots \\
   \end{pmatrix}.
\end{eqnarray}
We call $C$ and $\alpha$ the shift matrix and shift coefficient, respectively.
Applying row operations to the above matrix, it becomes the following matrix
\begin{eqnarray}
\label{eq:sec1_5}
   \begin{pmatrix}
        1  & -1 + 2\alpha & 1 & -1 + 2\alpha & \cdots  \\
        -1 + 2\alpha & 1 & -1 + 2\alpha & 1 & \cdots  \\
        0 & 0 & 0 & 0 &\cdots  \\
        0 & 0 & 0 & 0 &\cdots  \\
        \vdots & \vdots & \vdots & \vdots & \ddots \\
   \end{pmatrix}.    
\end{eqnarray}

A noteworthy feature is that the rank of $(1 - \alpha) B + \alpha C$ is discontinuous depending on $\alpha$: it is $1$ at $\alpha = 0$ and $\alpha = 1$, yet $2$ for every other value of $\alpha$.
We note that $\alpha$ values for rank $1$ matrices correspond to the solutions of the equation
\begin{eqnarray}
\label{eq:sec1_5}
   1 - (-1 + 2\alpha)^2 = 0.
\end{eqnarray}
When $\alpha = 0$, two distinct colors are present, whereas when $\alpha = 1$, only a single color remains. In this sense, $\alpha = 0$ represents a special case in which the rank is reduced with two colors.

\section{Matrix Decomposition}

Although matrices defined over $\{1,0\}$ and $\{1,-1\}$ may appear superficially identical, a matrix decomposition reveals that they differ significantly in their numerical behavior.

First, clearly, the rank of matrix $C$ is $1$.
Since the ranks of $B$ and $C$ are $1$,
$B$ and $C$ can be expressed as follows:
\begin{eqnarray}
\label{eq:sec2_3}
   B = \sigma \bm{u} \bm{v}^\mathsf{T}, \ 
   C = \tau \bm{w} \bm{z}^\mathsf{T}.
\end{eqnarray}
Here, $\bm{u}$ and $\bm{w}$ are $m$-dimensional vectors of norm 1, $\bm{v}^\mathsf{T}$ and $\bm{z}^\mathsf{T}$ are the transposes of $n$-dimensional vectors of norm 1, and $\sigma$ and $\tau$ are scalar values corresponding to $B$ and $C$, respectively.
Then, 
\begin{eqnarray}
\label{eq:sec2_4}
   A 
   = {1 \over 2} (B + C) 
   = {1 \over 2} (\sigma \bm{u} \bm{v}^\mathsf{T} + \tau \bm{w} \bm{z}^\mathsf{T}).
\end{eqnarray}
Since 
\begin{eqnarray} \nonumber
\label{eq:sec2_5}
   \bm{u}&=&({1 \over \sqrt{m}}, -{1 \over \sqrt{m}}, \cdots, {(-1)^{m-1} \over \sqrt{m}})^\mathsf{T}, \
   \bm{v}=({1 \over \sqrt{n}}, -{1 \over \sqrt{n}}, \cdots, {(-1)^{n-1} \over \sqrt{n}})^\mathsf{T}, \\ 
   \bm{w}&=&({1 \over \sqrt{m}}, {1 \over \sqrt{m}}, \cdots, {1 \over \sqrt{m}})^\mathsf{T}, \
   \bm{z}=({1 \over \sqrt{n}}, {1 \over \sqrt{n}}, \cdots, {1 \over \sqrt{n}})^\mathsf{T}, 
\end{eqnarray}
$\bm{u}$ and $\bm{w}$ are independent, and $\bm{v}$ and $\bm{z}$ are also independent. 
Thus, $A$ has rank of $2$.

Applying the singular value decomposition to $m \times n$ matrices $B$ and $C$ with equation (\ref{eq:intro_2}), we can obtain the singular value decomposition of $A$.
What above mentioned is consistent with the results of the singular value decomposition of $A$.
This can be illustrated by using an example of singular value decomposition for $4 \times 5$ matrix of
$A = U \Sigma V^\mathsf{T}$.
\begin{eqnarray}
\label{eq:sec2_1}
   A = U \Sigma V^\mathsf{T},
\end{eqnarray}
where, 
\begin{eqnarray}
\label{eq:sec2_2}
   A &=& 
   \begin{pmatrix}  
        1 & 0 & 1 & 0 & 1  \\
        0 & 1 & 0 & 1 & 0  \\
        1 & 0 & 1 & 0 & 1  \\
        0 & 1 & 0 & 1 & 0  \\
   \end{pmatrix}, \\
   U &=& 
   \begin{pmatrix}  
        {1 \over \sqrt{2}} & 0 & 0 & -{1 \over \sqrt{2}}  \\
         0 & {1 \over \sqrt{2}} & -{1 \over \sqrt{2}} & 0 \\
        {1 \over \sqrt{2}} & 0 & 0 & {1 \over \sqrt{2}}  \\
         0 & {1 \over \sqrt{2}} & {1 \over \sqrt{2}} & 0 \\
   \end{pmatrix}, \\
   \Sigma &=& 
   \begin{pmatrix}  
        \sqrt{6} & 0 & 0 & 0 & 0  \\
        0 & 2 & 0 & 0 & 0  \\
        0 & 0 & 0 & 0 &0  \\
        0 & 0 & 0 & 0 &0  \\
   \end{pmatrix}, \ \\
   V &=& 
   \begin{pmatrix}  
        {1 \over \sqrt{3}} & 0 & -{1 \over \sqrt{2}} & 0 & -{1 \over \sqrt{6}}  \\
         0 & {1 \over \sqrt{2}} & 0 & -{1 \over \sqrt{2}} & 0  \\
        {1 \over \sqrt{3}} & 0 & 0 & 0 & \sqrt{2 \over 3} \\
         0 & {1 \over \sqrt{2}} & 0 & {1 \over \sqrt{2}} & 0  \\
        {1 \over \sqrt{3}} & 0 & {1 \over \sqrt{2}} & 0 & -{1 \over \sqrt{6}} \\
   \end{pmatrix}.
\end{eqnarray}

\section{Information Embedding and Retrieving}

Although it is convenient to use the binary data $\{1,0\}$ when transmitting or storing digital data, it is more convenient to use $\{1,-1\}$ in computations because it makes it easier to determine things like orthogonality in linear algebra and correlation in statistics. Therefore, these code sets are often used in concert, exploiting the characteristics of each in their computational and communication roles.

Below is a simple example of embedding information within a checkered pattern and retrieving it.
We embed a small $2 \times 2$ matrix  into a $30 \times 30$ checkered pattern matrix, 
and all the entries of the small matrix is the mean value of checkered pattern matrix. 
If we use $\{1, 0\}$ for the checkered pattern, the mean value is $\displaystyle 0.5$, while using $\{1, -1\}$, it is $0$.

\begin{figure}[htb]
\begin{center}
\includegraphics[scale=0.4]{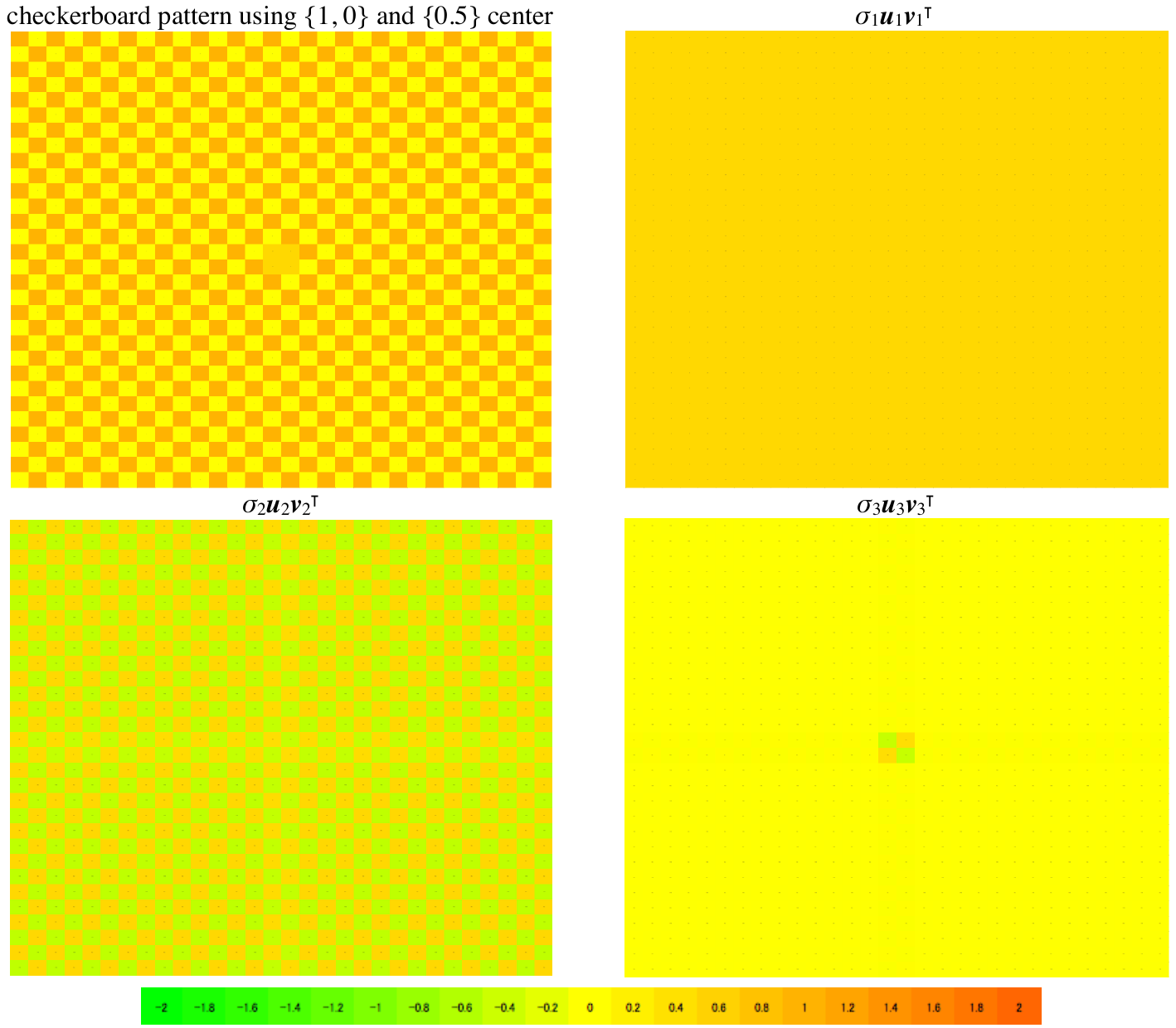}
\end{center}
\caption{Checkered pattern using $\{1, 0\}$ and $\{0.5\}$ in the center.} 
\label{fig:plus1zero}
\end{figure}

\begin{figure}[htb]
\begin{center}
\includegraphics[scale=0.4]{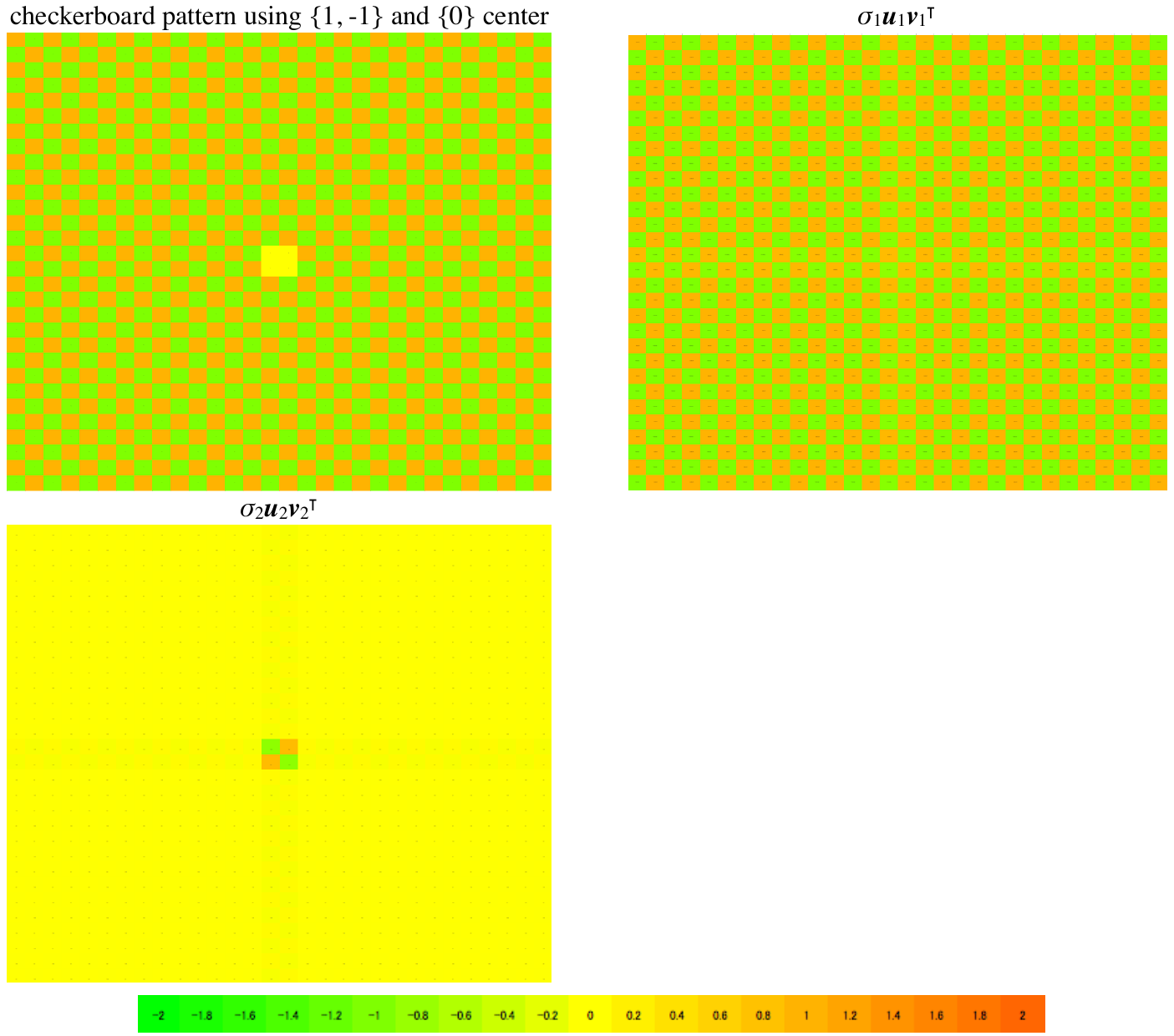}
\end{center}
\caption{Checkered pattern using $\{1, -1\}$ and $\{0\}$ in the center.} 
\label{fig:plus1minus1}
\end{figure}

Figure \ref{fig:plus1zero} shows the embedding of a small uniform matrix into a $\{1, 0\}$ checkered pattern. The original pattern is shown in the top left of the figure. The rest of the figure shows the matrices $\sigma_i \bm{u}_i \bm{v}_i^\mathsf{T} \ (i=1,2,3)$ generated using singular value decomposition.
Figure \ref{fig:plus1minus1} shows the embedding of a small uniform matrix into a $\{1, -1\}$ checkered pattern. The original pattern is shown in the top left of the figure. The rest of the figure shows the matrices $\tau_i \bm{w}_i \bm{z}_i^\mathsf{T} \ (i=1,2)$ generated using singular value decomposition.
Comparing Figures \ref{fig:plus1zero} and \ref{fig:plus1minus1}, we can see that in Figure \ref{fig:plus1zero}, the checkered mean of $\displaystyle 0.5$ is added to the entire matrix, which causes an increase in rank by 1 compared to the case in Figure \ref{fig:plus1minus1}.
In both cases, using singular value decomposition, the embedded small matrix is separated from the checkered pattern matrix, revealing the added pattern.

This example shows that it is more efficient to use $\{1, -1\}$ in matrix computations, but more convenient to use $\{1, 0\}$ in digital manipulations because the image data has non-negative values.

\section{Conclusion}

Checkered patterns are characterized by their square structure and the use of only two distinct colors. These colors are typically represented by two types of numerical sets: $\{1,0\}$ and $\{1,-1\}$. Although matrices based on $\{1,0\}$ may seem identical to those based on $\{1,-1\}$ when forming checkered patterns, these two kinds of matrices are completely different in a mathematical sense because a matrix using $\{1, 0\}$ has a rank of $2$ and a matrix using $\{1, -1\}$ has a rank of $1$. The difference in rank $1$ occurs because the mean value of the two numbers is either non-zero or zero. However, this property of rank reduction can be exploited in matrix operations to decrease computational effort such as information embedding and retrieving.




\end{document}